\documentclass[graybox]{svmult}

\usepackage{type1cm}        
\usepackage{makeidx}         
\usepackage{graphicx}        
\usepackage{multicol}        
\usepackage[bottom]{footmisc}

\usepackage{newtxtext}       %
\usepackage{newtxmath}       

\makeindex             

\usepackage{multirow}
\usepackage{soul}
\usepackage{tabularx,ragged2e}
\newcolumntype{C}{>{\Centering\arraybackslash}X} 
\usepackage{caption}
\usepackage{subcaption}
\usepackage{float}
\usepackage{xcolor}
\usepackage{pdfpages}

\newcommand{\bsu}{\mathsf{\mathbf{u}}}
\newcommand{\bsv}{\mathsf{\mathbf{v}}}

\newcommand{\bsf}{\mathsf{\mathbf{f}}}

\newcommand{\beq}{\begin{equation}}
\newcommand{\eeq}{\end{equation}}

\newcommand{\cA}{{\mathcal A}}

\newcommand{\cI}{{\mathcal I}}

\newcommand{\cQ}{\mathcal{Q}}

\newcommand{\cR}{\mathcal{R}}  
\newcommand{\tcR}{\widetilde{\mathcal{R}}}

\newcommand{\cT}{{\mathcal T}}

\newcommand{\cV}{{\mathcal V}}

\newcommand{\cB}{{\mathcal B}}

\newcommand{\bx}{\boldsymbol{x}}

\newcommand{\supp}{\mathrm{supp}}

\newcommand{\ri}{{\rm i}}

\usepackage{color}
\usepackage{hyperref}
\definecolor{myblue}{rgb}{0,0,0.6}
\definecolor{darkgreen}{rgb}{0,0.5,0}
\hypersetup{colorlinks=true,
linkcolor=myblue,citecolor=myblue,filecolor=myblue,urlcolor=myblue}

\definecolor{escol}{rgb}{0,0.4,0}
\definecolor{sgcol}{rgb}{0,0,0.7}
\definecolor{estcol}{rgb}{0.5,0,0}
\definecolor{esnewcol}{rgb}{0,0.5,0}

 \newcommand{\igg}[1]{{\color{black}{#1}}}

\newcommand{\beqs}{\begin{equation*}}
\newcommand{\eeqs}{\end{equation*}}
\newcommand{\bit}{\begin{itemize}}
\newcommand{\eit}{\end{itemize}}
\newcommand{\ben}{\begin{enumerate}}
\newcommand{\een}{\end{enumerate}}
\newcommand{\bal}{\begin{align}}
\newcommand{\eal}{\end{align}}
\newcommand{\bals}{\begin{align*}}
\newcommand{\eals}{\end{align*}}
\newcommand{\bre}{\begin{remark}}
\newcommand{\ere}{\end{remark}}
\newcommand{\bpf}{\begin{proof}}
\newcommand{\epf}{\end{proof}}
\newcommand{\ble}{\begin{lemma}}
\newcommand{\ele}{\end{lemma}}
\newcommand{\bco}{\begin{corollary}}
\newcommand{\eco}{\end{corollary}}
\newcommand{\bex}{\begin{example}}
\newcommand{\eex}{\end{example}}
\newcommand{\bth}{\begin{theorem}}
\newcommand{\enth}{\end{theorem}}

\newcommand{\mymatrix}[1]{\mathsf{#1}}
\newcommand{\MA}{{\mymatrix{A}}}
\newcommand{\MB}{{\mymatrix{B}}}
\newcommand{\MD}{{\mymatrix{D}}}
\newcommand{\MR}{{\mymatrix{R}}}
\newcommand{\bMR}{{\mymatrix{\mathbf{R}}}}
\newcommand{\bMQ}{{\mymatrix{\mathbf{Q}}}}
\newcommand{\bMT}{{\mymatrix{\mathbf{T}}}}

\newcommand{\MI}{{\mymatrix{I}}}
\newcommand{\ME}{{\mymatrix{E}}}

\newcommand{\tfa}{\text{ for all }}

\newcommand{\tand}{\text{ and }}

\begin{document}

\title*{A variational interpretation of Restricted Additive Schwarz with impedance transmission condition  for the Helmholtz problem}
\titlerunning{RAS with local impedance solves for Helmholtz}
\author{Shihua Gong, Martin J.~Gander,  Ivan G.~Graham and Euan A.~Spence}
\authorrunning{S. Gong, M. J.~Gander, I. G.~Graham,  and E. A.~Spence}
\institute{Shihua Gong, Ivan G.~Graham and Euan A.~Spence \at  Department of Mathematical Sciences, University of Bath, Bath BA2
7AY, UK.
\and Martin J.~Gander \at Department of Mathematics, University of Geneva, Switzerland.}
%
%
\maketitle

\abstract{ In this paper we revisit the Restricted Additive Schwarz method for solving discretized Helmholtz problems, using  impedance boundary conditions on subdomains (sometimes called ORAS). We present this method in its variational form and show that it   can be seen  as a finite element discretization of a  parallel overlapping domain decomposition method defined at the PDE level.
  In a recent companion paper, the authors have  proved certain contractive properties of the error propagation operator for this method at the PDE level, under certain geometrical assumptions. We illustrate computationally that these properties are also enjoyed by its finite element approximation, i.e.,  the ORAS method. }

\section{The Helmholtz problem} 
Motivated by  the large range of applications,  there is currently great interest in designing and analysing preconditioners for finite element discretisations of the Helmholtz equation
\begin{align} \label{eq:Helm}
 - (\Delta  + k^2 \, )u  \ = \ f  \quad \text{on} \quad \Omega\, ,  
  \end{align}  
  on a $d-$dimensional domain $\Omega$  ($d = 2,3$),  with $k$ the ({assumed constant, but} possibly large) angular frequency.   While the methods presented  easily apply to quite general scattering problems and geometries,  we restrict attention here  to the interior impedance problem,  where $\Omega$ is bounded,  and the boundary condition is
  \begin{align}
\left(\frac{    \partial } {\partial n} - \ri k\right) u = g \quad \text{on} \quad
     \partial \Omega\, , \label{eq:BC}
  \end{align}
where $\partial u/\partial n$ is the outward-pointing normal derivative of $u$ on  $\Omega$.

    The weak form of problem \eqref{eq:Helm}, \eqref{eq:BC} is to seek $u \in H^1(\Omega)$ such that
\begin{align} \label{eq:weak} a(u,v) = F(v) : =   \int_{\Omega} f\bar{v}~dx + \int_{\partial \Omega}g \bar{v} ~ds,
\end{align} 
     \begin{equation*}
  \text{where} \quad \quad     a(u,v) := \int_\Omega(\nabla u . \nabla \overline{v}  - k^2 u \overline{v}) - \ri k \int_{\partial \Omega}  u \overline{v} , \quad \text{for}\quad u,v \in H^1(\Omega) .
    \end{equation*}

 \section{Parallel iterative Schwarz method}

To solve  \eqref{eq:Helm}, \eqref{eq:BC},  we  shall consider domain decomposition methods,    based on a set of Lipschitz polyhedral    
subdomains   $\{\Omega_\ell\}_{\ell = 1}^N$,  
forming  an overlapping cover of $\Omega$ and equipped with  a partition of unity:
 $\{\chi_\ell\}_{\ell = 1}^N$, such that 
 \beq
 \left.
 \begin{array}{ll} 
  \text{for each} \quad \ell: & \  \supp \chi_\ell \subset \overline{\Omega_\ell}, \quad
                                0 \leq  \chi_\ell(\bx) \leq 1\,\,  \text{when } \bx \in \overline{\Omega_\ell},
    \\
    & \\
    \quad\tand\quad  & 
   \sum_\ell  \chi_\ell(\bx) = 1 \,\tfa \bx \in \overline{\Omega}.
  \end{array}
  \right\}
 \label{POUstar}
  \eeq

Then,   the parallel Schwarz method for \eqref{eq:Helm}, \eqref{eq:BC}  {with Robin (impedance) transmission conditions  is: given  $u^n$ defined on $\Omega$,  we solve the local problems:}  
\begin{align}
  -(\Delta + k^2)u_\ell^{n+1}  & =  f  \quad  &\text{in } \  \Omega_\ell\, ,  \label{eq21}\\
  \left(\frac{\partial }{\partial n_\ell} - \ri k \right) u_\ell^{n+1}  & = \left(\frac{\partial }{\partial n_\ell} - \ri k \right) u^n \quad  &\text{on }  \ \partial \Omega_\ell\backslash \partial \Omega\, , \label{eq22}  \\
    \left(\frac{\partial }{\partial n_\ell} - \ri k \right) u_\ell^{n+1}  & = g  \quad  &\text{on } \  \partial \Omega_\ell  \cap \partial \Omega.   \label{eq23} 
\end{align}
Then the next  iterate is the weighted sum of the local solutions
\begin{align} \label{star}
  u^{n+1} := \sum_\ell \chi_\ell  u_\ell ^{n+1}.  
\end{align}
{Information is shared between neighbouring subdomains at each iteration via 
\eqref{star}.}

{In   \cite{GoGaGrLaSp:21}}, we analyse the iteration \eqref{eq21} -- \eqref{star} in the
function space 
   $$ U (\Omega) := \big\{ v \in H^1(\Omega): \,  \Delta v  \in L^2(\Omega), \, \partial v / \partial n  
   \in L^2(\partial \Omega) \big\},   $$
   {and  its local analogues   $U(\Omega_{\ell})$.} 
   Using the fact that any function  $ v \in U(\Omega_\ell)$ has  impedance trace $(\partial/\partial n  - \ri k)v \in L^2(\Gamma)$ on any Lipschitz curve $\Gamma \subset \Omega_\ell$, we prove in \cite{GoGaGrLaSp:21}  that  \eqref{eq21} -- \eqref{star} is well-defined in the space $U(\Omega)$. {Moreover, {introducing}  $e_\ell^n = u\vert_{\Omega_\ell} - u^n_\ell$,   and letting  $\mathbf{e}^n = (e_1^n, \ldots, e_N^n)$, we prove in  \cite{GoGaGrLaSp:21} that $\mathbf{e}^{n+1}  = \cT \mathbf{e}^n, $  where  under certain geometric assumptions,   $\cT$ has the `power contraction' property    \begin{align}\label{eq:errorprop}  \Vert \cT^N \Vert \ll 1, \end{align} with respect to the product norm on   $\prod_\ell U_0(\Omega_\ell)$, where $U_0(\Omega_\ell)$ is the subspace of functions  $v \in U(\Omega_\ell)$,  for which $\Delta v + k^2 v = 0$ on $\Omega_\ell$.  \igg{Analogously to  \cite{BeDe:97}, the  norm of $v$ is the $L^2$ norm of its impedance data on $\partial \Omega_\ell$.}  See the remarks in \S \ref{sec:Numerical}, especially \eqref{powerbdd},    for a more precise explanation of \eqref{eq:errorprop}.


    
   The aim of this note is to show that a natural finite element analogue of  \eqref{eq21} -- \eqref{star} corresponds to a preconditioned Richardson-type iterative method for the finite element approximation of  \eqref{eq:Helm}, \eqref{eq:BC}, where the preconditioner is  a {Helmholtz-orientated} version of the popular  Restricted Additive Schwarz method. 
This preconditioner is  given several different names in the literature  --  {WRAS-H} (Weighted RAS for Helmholtz) \cite{KiSa:07}, ORAS (Optimized Restricted Additive Schwarz) \cite{st2007optimized,DoJoNa:15,GoGrSp:20},  IMPRAS1 (RAS with impedance boundary condition) \cite{GrSpVa:17a}. {However it has not previously been directly connected {via a variational argument}  to the iterative method  \eqref{eq21} -- \eqref{star} in the Helmholtz case, {although there are algebraic discussions  (e.g., \cite{EfGa:03}, \cite[\S2.3.2]{DoJoNa:15})}. We also demonstrate numerically in \S \ref{sec:Numerical}, that the finite element analogue of  \eqref{eq21} -- \eqref{star} inherits the property \eqref{eq:errorprop} proved at the continuous level in  \cite{GoGaGrLaSp:21}. }

\igg{Method  \eqref{eq21}--\eqref{star} is an example of methods 
  studied more generally in the Optimized Schwarz literature (e.g., \cite{gander2006optimized,st2007optimized}), 
  where Robin  (or more sophisticated) transmission  conditions
  are constructed with the aim of  optimizing  convergence rates.
  Although  the transmission condition \eqref{eq22} above can be justified directly as a
  first order absorbing condition for the local Helmholtz problem \eqref{eq21} (without considering optimization), this method 
  is still often called `Optimized Restricted Additive Schwarz' (or `ORAS') and we shall continue  this naming  convention here.}
ORAS is arguably the most successful one-level parallel method for Helmholtz problems. It can be applied on very general geometries, does not depend on parameters, and can even be robust to increasing $k$ \cite{GoGrSp:20}. More generally it can be combined with coarse spaces to improve its robustness properties. }

\section{Variational formulation of  RAS with impedance transmission condition (ORAS)  }
{Here we formulate a  finite element approximation of \eqref{eq:Helm}, \eqref{eq:BC} and show that it coincides with ORAS.  We  introduce 
a  nodal finite element space $\cV^h \subset H^1(\Omega)$ consisting of continuous piecewise  polynomials of total degree  $\leq p$ on a conforming mesh $\cT^h$. 
 Functions in  $\cV^h$  are  uniquely determined by their values at
 nodes in $\overline{\Omega}$,  denoted $\{ x_j: j \in \cI\}$,  for some index set $\cI$.  
The local space on $\overline{\Omega_\ell}$ is    
 $\cV_{\ell}^h: =\{v_{h}|_{\overline{\Omega_\ell}}~:~v_{h}\in {\cV}^h\}$ with corresponding nodes  denoted $\{x_j: j \in \cI_\ell\}$, for some  $\cI_\ell \subset \cI$.   }

 Using the sesquilinear form $a$ and right-hand side $F$ appearing in \eqref{eq:weak},   we can define the discrete
 operators $\cA_h, F_h :\cV^h \mapsto (\cV^h)'$  
 by
\begin{equation}
  (\cA_h u_h)(v_h) : = a(u_h,v_h)\quad \text{and}   \quad   F_h(v_h) =  F(v_h), \quad \text{for all} \quad u_h , v_h \in \cV_h .\label{discrete}  
  \end{equation}
Analogously,  on each subdomain $\Omega_\ell$,  we define $\cA_{h,\ell}: \cV_\ell^h \rightarrow (\cV_{\ell}^h)'$ by
$(\cA_{h,\ell}u_{h,\ell})(v_{h,\ell}) : = a_\ell(u_{h,\ell}, v_{h,\ell}).$
We also need prolongations $\cR_{h,\ell} ^\top, \tcR_{h, \ell}^\top : \cV_\ell^h \rightarrow \cV^h$ defined for all $v_{h,\ell} \in \cV^h_\ell$ by
\begin{equation}\label{nodewise}
(\cR_{h,\ell} ^\top  v_{h,\ell})(x_j)  = \left\{ \begin{array}{ll} v_{h,\ell} (x_j) & \quad {j \in \cI_\ell}, \\ 
                                          0 & \quad \text{otherwise, }\end{array} \right.  \quad \text{and} \quad \tcR_{h,\ell}^\top v_{h,\ell}  = \cR_{h,\ell}^\top (\chi_\ell v_{h,\ell}). 
\end{equation}

{Note the subtlety in \eqref{nodewise}: The extension  $\cR_{h,\ell}^\top v_{h,\ell}$ is defined {\em nodewise}: It  coincides with $v_{h,\ell}$ at  nodes in  $\overline{\Omega_\ell}$ and vanishes at nodes in $\Omega\backslash \overline{\Omega_\ell}$. Thus $\cR_{h,\ell} ^\top  v_{h,\ell} \in \cV^h \subset H^1(\Omega)$. This is  \igg{an $H^1-$}  conforming finite element approximation of the zero extension of $v_{h,\ell}$ to all of ${\Omega}$.  (The zero extension is  not in $H^1(\Omega)$ in general.) }  
We define the restriction operator $\cR_{h,\ell}: \cV_h' \rightarrow \cV_{h, \ell}'$  by duality, i.e., for all
$F_h \in \cV_h'$,  $$ (\cR_{h,\ell}F_h)(v_{h,\ell}) := F_{h}(\cR_{h,\ell}^\top v_{h,\ell}), \quad v_{h, \ell} \in \cV_\ell^h .$$ 
                                      Then the ORAS preconditioner is  the operator $\cB_h^{-1}  : \cV_h' \rightarrow \cV_h$ defined by
\begin{align} \label{eq:PC} \cB_h^{-1} \ : =\  \sum_\ell \tcR_{h,\ell}^\top \cA_{h,\ell}^{-1}\cR_{h,\ell}. 
\end{align}
{This preconditioner can also be written in terms of operators $\cQ_{h, \ell} : \cV^h \rightarrow \cV_\ell^h$ defined for all $u_h \in \cV^h$  by} 
\begin{align} a_\ell(\cQ_{h, \ell} u_h , v_{h, \ell} ) \ = \ a (u_h, \cR_{h,\ell}^\top v_{h,\ell}), \quad
  \text{for all} \quad v_{h,\ell} \in \cV^h_\ell ,  \label{eq:var}
\end{align}
where   $\cR^\top_{h,\ell}$ is defined in \eqref{nodewise},  and then 
 $\cB_h^{-1} = \sum_\ell \tcR^{\top}_{h,\ell} \cQ_{h,\ell}$. 
The corresponding preconditioned Richardson iterative method  can be written as
     \begin{equation}\label{eq:oras-iter}
     u^{n+1}_h = u^n_h  + \cB_h^{-1}  (F_h - \cA_h u_h^n).
     \end{equation}
 The matrix realisation of \eqref{eq:oras-iter} is given in  \S \ref{sec:Numerical}.


 \section{Connecting  the parallel iterative method with ORAS}

In this section, we show that a natural finite element approximation of \eqref{eq21}--\eqref{star} yields \eqref{eq:oras-iter}. 
First, to write  \eqref{eq21} - \eqref{star} in  a residual correction form, we   introduce the   ``corrections''  $\delta_\ell^n : = u_{\ell}^{n+1} - u^n\vert_{\overline{\Omega_\ell}}, $. With this definition we have   
\begin{align}
  -(\Delta + k^2)\delta_\ell^{n}  & =  f +   (\Delta  + k^2) u^{n} \quad  \text{in } \  \Omega_\ell\, ,  \label{eq21c}\\
  \left(\frac{\partial }{\partial n_\ell} - \ri k \right)\delta_\ell^{n}  & =0  \quad  \text{on }  \ \partial \Omega_\ell\backslash \partial \Omega\, , \label{eq22c}  \\
                                                                                                                                                                                     \left(\frac{\partial }{\partial n_\ell} - \ri k \right)\delta_\ell^{n}  & = g -   \left(\frac{\partial }{\partial n_\ell} - \ri k \right)  u^{n} \quad  \text{on } \  \partial \Omega_\ell  \cap \partial \Omega\, ,  \label{eq23c}                                                          
\end{align}
\begin{align}
  \text{and then  } \quad \quad & u^{n+1} = u^{n} + \sum_{\ell} \chi_\ell \delta_\ell^{n}.    \label{eq:sum-correction}
\end{align} 
Note, there is more  subtlety here: Because of \eqref{star},  $u^n\vert_{\overline{\Omega_\ell}}$ is {\em not} the same as $u^n_\ell$. 
The theory in \cite{GoGaGrLaSp:21} {can be used to show that \eqref{eq21c}--\eqref{eq:sum-correction} is   still}  well-posed in $U(\Omega)$.              Multiplying \eqref{eq21c} by $v_\ell\in H^1(\Omega_\ell)$, integrating by parts and using  \eqref{eq22c}, \eqref{eq23c}, $\delta_\ell^n$ satisfies,  for $ v_\ell \in H^1(\Omega_\ell)$,

          \begin{align}\nonumber 
            a_\ell(\delta_\ell^n , v_\ell) & \ =\   \int_{\Omega_\ell}f\, \overline{ v_\ell}  \  +\  \int_{\partial \Omega_\ell \cap \partial \Omega} g
                                             \, \overline{v_\ell}  \\
                                           & \ + \ \ \int_{\Omega_\ell}(\Delta + k^2) u^n\,  \overline{ v_\ell} \  -\ 
                                             \int_{\partial \Omega_\ell \cap \partial \Omega}
                                 \left( \frac{\partial }{ \partial n_\ell} - \ri k \right) u^n \, \overline{ v_\ell} .  \label{eq:local-correction-weak1}
          \end{align}
         { To implement the   finite element discretization of this,  we will need to handle the case when     $u^n$ on the right-hand side is replaced by a given iterate $u_h^n \in \cV^h$ and when   
the  test function
          \igg{$v_\ell \in H^1(\Omega_\ell)$} is replaced by  $v_{h,\ell} \in \cV^h_\ell$.   } The  third term on the right hand side of \eqref{eq:local-correction-weak1} then requires integration by parts to make sense. {Using the  nodewise extension $\cR^\top_{h, \ell}$  we replace the third and fouth terms in     \eqref{eq:local-correction-weak1} by}
          \begin{align} \label{approx}
\int_{\Omega}(\Delta + k^2) u^n\,  \overline{ \cR_{h,\ell}^\top  v_{h,\ell}} \  -\ 
                                             \int_{\partial \Omega}
                                             \left( \frac{\partial }{ \partial n} - \ri k \right) u^n \,\,  \overline{ \cR_{h,\ell}^\top
                                               v_{h,\ell}} \  = \ -a(u^n,\cR_{h,\ell}^\top
                                             v_{h,\ell}),   \end{align}
                                             where the right-hand side  is obtained from the left via  integration by parts over $\Omega$.                    This leads to the  FEM analogue of \eqref{eq21c} -- \eqref{eq:sum-correction}:        Suppose $u_h^n \in \cV^h$ is given. Then 
                                             \begin{align}\label{twenty}u_h^{n+1} : = u_h^n + \sum_\ell \tcR_{h,\ell}^\top \delta_h^n, \end{align} where (using   
                                             \eqref{eq:local-correction-weak1},  \eqref{approx} and \eqref{discrete}),
                                             $                                             a_\ell(\delta_{h,\ell}^n, v_{h,\ell}) \ = \ \cR_{\ell,h} F_h \ -\
                                             a(u^n,\cR_{\ell,h}^\top
                                             v_{\ell,h}).$
Thus,
$$ \delta_{h,\ell}^n \ = \ \cA_{h,\ell}^{-1} \cR_{h,\ell} (F_h - \cA_h u_h^n).$$
Combining this with    \eqref{twenty}, we obtain exactly  \eqref{eq:oras-iter}.

     \section{Numerical results}
     \label{sec:Numerical}
     
     Denoting the nodal bases for $\cV^h$ and $\cV_\ell^h$ by 
     $\{\phi_j\}$
     and  $\{\phi_{\ell, j}\} $ respectively, we introduce stiffness matrices
  $
     \MA_{i,j} : =  a(\phi_j,\phi_i)$  \text{and} $(\MA_\ell)_{i,j} :  = a_\ell(\phi_{\ell,j}, \phi_{\ell,i})$,    and the load vector $f_i := F_h(\phi_i)$. Then  we can write  \eqref{eq:oras-iter} as
  \begin{equation}\label{eq:oras-matrix}
     \bsu^{n+1} = \bsu^n + \MB^{-1}  (\bsf - \MA \bsu^n).
   \end{equation}
   Here 
   $\bsu^n$ is the coefficient vector of $u_h^n$ with respect to the nodal basis of $\cV^h$,  and 
      $$\MB^{-1} = \sum_{\ell}\tilde{\MR}_\ell^\top \MA_\ell^{-1}\MR_\ell\,  , $$
where 
     $
     (\MR_\ell^\top)_{p,q} := (\cR_\ell^\top \phi_{\ell,q})(x_p)$,  $(\tilde{\MR}_\ell^\top)_{p,q} := (\tcR_\ell^\top \phi_{\ell,q})(x_{\ell,p}),$
 and $\MR_\ell =  (\MR_\ell^\top)^\top.$

    

     In this section, ({motivated by \eqref{eq:errorprop}}),  we numerically investigate the contractive property of the ORAS iteration \eqref{eq:oras-matrix}. Letting $\bsu$ be the solution of 
     $
     \MA\bsu = \bsf,
     $
we can combine with \eqref{eq:oras-matrix} to  obtain the error propagation equation 
$$
\bsu^{n+1} - \bsu = \ME(\bsu^n -\bsu), \quad \text{where} \quad  \ME= \MI - \MB^{-1}A .
$$ 
{Since $\sum_\ell \tilde{\MR}^\top_\ell \MR_\ell = \MI$,    we can write }
$$\ME= \sum_{\ell} \tilde{\MR}_\ell^\top (\MR_\ell - \MA_\ell^{-1} \MR_\ell \MA) = \tilde{\bMR}^\top (\bMR - \bMQ), $$
 where $\tilde{\bMR}^\top$ is the row  vector of  matrices:  {$\tilde{\bMR}^\top = (\tilde{\MR}_1^{\top},\tilde{\MR}_2^{\top},\cdots, \tilde{\MR}_N^{\top})$, and $\bMR = (\MR_1; \MR_2;\cdots;\MR_N)$, and $\bMQ = (\MA_1^{-1} \MR_1 \MA; \MA_2^{-1} \MR_2 \MA,\cdots, \MA_N^{-1} \MR_N \MA )$ are column vectors. 
 Then it is easily seen that   $\ME\tilde{\bMR}^\top = \tilde{\bMR}^\top \bMT$,  where  $ \bMT :=  (\bMR - \bMQ) \tilde{\bMR}^\top$. {Moreover, since $\tilde{\bMR}^\top \bMR = \MI$, we have $\bMT \bMR \tilde{\bMR}^\top = \bMT$,  and so     it follows that}
 \begin{align} \label{powerE}
\ME^s = \tilde{\bMR}^\top \bMT^{s} \bMR \quad \text{for any} \quad s\geq 1, \quad  
 \end{align}
 {As explained in \cite[\S5.1]{GoGaGrLaSp:21},}    $\bMT$ is a discrete  version of the operator  $\cT$ appearing in \eqref{eq:errorprop} above. \igg{In \cite{GoGaGrLaSp:21}, we study  fixed point iterations with matrix
   $\bMT$ and use these   to illustrate various properties of   the fixed point  operator $\cT$ 
   in the product norm described above. In this paper we consider only the norms of $\ME^s$.  By \eqref{powerE},  if  $\bMT^s$ is sufficiently contactive,   then $\ME^s$ will also be contractive.}   
 
 To compute the norm of $\ME^s$, we introduce the vector norm:   
$
\|\bsu \|_{1,k}^2 =  \bsu^* \MD_{k} \bsu 
$, for    $\bsu\in\mathbb{C}^M,$ 
where  $M = \mathrm{dim}(\cV^h)$ and,   for all nodes $x_p, x_q$ of $\cV^h$,    
$(\MD_{k})_{p,q} = \int_{\Omega} \nabla \phi_p \cdot \nabla \phi_q + k^2\phi_p\phi_q~dx,$ . This is the matrix induced by the usual $k-$weighted $H^1$ inner product on $\cV^h$ .  
We shall  compute 
$$
\|\ME^s\| : =\max_{0\neq\bsv\in\mathbb{C}^M }\frac{ \| \ME^s \bsv\|_{1,k}}{\|\bsv\|_{1,k}}, \quad \text{for integers} \quad s \geq 1, 
$$
{which is equal to the square root of the largest eigenvalue of the matrix $\MD_{k}^{-1}(\ME^*)^s \MD_{k}\ME^s$. This is computed using the SLEPc facility within the package FreeFEM++ \cite{hecht2019freefem++}.
In the following numerical experiments, done on rectangular domains, we use conforming Lagrange elements of degree  $2$, on  uniform meshes with mesh size decreasing  with $h\sim k^{-5/4}$ as $k$ increases, sufficient for avoiding the pollution effect. }

{We consider two different examples of domain decomposition.  First we consider a long  rectangle of size $(0,\frac{2}{3}N )\times (0,1)$, partitioned  into $N$ non-overlapping strips of equal width $2/3$. We  then extend each subdomain  by adding neighbouring elements whose distance from the boundary is $\leq {1}/{6}$.  
  This gives an  overlapping  cover, with each subdomain  a unit square, except for the  subdomains at the ends,  which are rectangles with aspect ratio $6/5$. }
For this example, a rigorous estimate ensuring \eqref{eq:errorprop} is proved in \cite{GoGaGrLaSp:21}. The result implies  that 
  \begin{align} \label{powerbdd}  \Vert \cT^N \Vert \ \leq C (N-1) \rho \ + \ \mathcal{O}(\rho^2). \ \end{align}
  Here,  $\rho$ is the maximum of the $L^2$ norms of the `impedance maps' which describe the exchange of impedance data between boundaries of overlapping subdomains within a single iteration. The constant $C$ is independent of $N$, but the hidden constant may depend on $N$. Thus for small enough $\rho$,  $\cT^N$ is a contraction.
  Conditions ensuring this are explored in \cite{GoGaGrLaSp:21}.
\begin{table}[t]
\centering
\scalebox{0.8}{
\setlength\extrarowheight{2pt} 
\begin{tabularx}{1.2\textwidth}{C|CCC|CCC|CCC|CCC}
\hline
\hline
$ N $& \multicolumn{3}{c}{$2$}& \multicolumn{3}{c}{$4$}& \multicolumn{3}{c}{$8$}& \multicolumn{3}{c}{$16$} \\ \cline{2-13}
k	& $\|\ME\|$ & $\|\ME^s\|$ & $\| \ME^{s+1}\|$	& $\|\ME\|$ & $\|\ME^{s-1}\|$	& $\|\ME^{s}\|$	& $\|\ME\|$ & $\|\ME^{s-1}\|$	& $\|\ME^{s}\|$		& $\|\ME\|$ & $\|\ME^{s-1}\|$	& $\|\ME^{s}\|$\\ 
\hline 
$20$		&5.6 &0.52 & 0.05 &5.8 &5.24 &0.18 &5.8 &4.5 & 0.11	&5.9	&3.4	&0.17\\
$40$		&9.0	&1.0 & 0.094 	&9.1	&8.5	&0.46	&9.1	&8.1	&0.34	&9.1	&7.6	&0.36\\
$80$		&14.3 &1.9 & 0.17 & 14.3 &13.1 &0.78	&14.3	&13.0&0.61	&14.3&12.6	&0.66   \\
\hline
	\hline
\end{tabularx}}
\caption{Strip partition of   $(0,\frac{2}{3}N )\times (0,1)$: Norms of powers of  $\ME$
  ($s = N$)  }\label{tb:norm-strip}
\end{table}

In Table \ref{tb:norm-strip} we observe  the rapid drop in the norm of $\Vert \ME^s \Vert $ compared with $\Vert \ME^{s-1}\Vert$ (with $s = N$). Moreover $\ME^N$ is a contraction when $N=4,8,16$. When  $N=2$ we do not have $\ME^2$ contracting, but $\ME^3$ certainly is. 
 Although  $\Vert E^N \Vert $ is increasing (apparently linearly) with $k$,
 $\Vert E^s\Vert$   decreases rapidly for $s> N$,  when  $k$ is fixed. 
 Note that  $\Vert \ME\Vert $ can be quite large, and is growing as $k$ increases:  thus the error of the iterative method may  grow initially before converging to zero. 
 \igg{
Also, although the right-hand side of \eqref{powerbdd} grows linearly in $N$ for fixed $\rho$, the norm of $E^N$ does not exhibit substantial growth. Thus we conclude that}  \eqref{powerbdd} may be pessimistic in its $N$-dependence. In fact sharper  estimates are proved and explored computationally in \cite{GoGaGrLaSp:21}. \igg{An interesting open question is to find a lower bound for $s$ as a function of $N$ and $k$ which ensures contractivity.}

   In \cite{GoGaGrLaSp:21} it is shown that the computation of  $\rho$, or related more detailed quantities can be done by solving eigenvalue problems on subdomains.  This, combined with estimates like   \eqref{powerbdd}   could be seen as an {\em a priori} condition for convergence,  rather like convergence predictions  via  condition number estimates. These  always give a sufficient  condition for good performance (which is often not sharp).  


{
In  the next experiment the domain $\Omega$ is the unit square,  divided into $N \times N $ equal square subdomains  in a ``checkerboard'' domain decomposition. Each subdomain is extended by adding neighbouring elements a distance $\leq 1/4$ of the width of the non-overlapping subdomains,  thus yielding an overlapping domain decomposition with ``generous'' overlap.
\begin{table}[h]
\centering
\scalebox{0.8}{
\setlength\extrarowheight{2pt} 
\begin{tabularx}{1.2\textwidth}{C|CC|CC|CC|CCC}
\hline
\hline
$ N\times N$& \multicolumn{2}{c}{$2\times 2$}& \multicolumn{2}{c}{$4\times 4$}& \multicolumn{2}{c}{$6\times 6$}& \multicolumn{3}{c}{$8\times 8$}  \\ \cline{2-10}
	k& $\|E^{s-1}\|$ & $\|E^{s}\|$	& $\|E^{s-1}\|$	& $\|E^{s}\|$	& $\|E^{s-1}\|$	& $\|E^{s}\|$		& $\|E^{s-1}\|$	& $\|E^{s}\|$& {\tt GMRES} \\ 
\hline 
$20$		&4.0e-1 &8.8e-2 &2.3e-3 &1.2e-3 &38 &41 &1.2e6 & 1.4e6	&34\\
$40$		&7.2e-1 &1.6e-1 &4.4e-2 &2.8e-3 &1.5e-3 &1.0e-3 & 6.4e-5 &5.3e-5& 28 \\
$80$		&1.0 &2.4e-1 & 1.5e-2 &9.8e-3 &3.9e-4&2.8e-4& {1.9e-6}&{9.2e-7}& 26  \\
$160$&	1.8	&5.0e-1	&1.1e-2	&6.3e-3	&7.3e-4	&5.3e-4 &{ 9.2e-5}& {7.5e-5}&  24  \\
\hline
	\hline
\end{tabularx}}
\caption{Checkerboard partition of the unit square:  Norms of  powers of $\ME$  ($s=N^2$), }\label{tb:norm-square1}
\end{table}
In Table \ref{tb:norm-square1}  we tabulate  $\Vert \ME^{s-1}\Vert$ and $\Vert \ME^s\Vert$,  for  $s = N^2$ (i.e., the total number of subdomains).  Here we do not see such a difference between these two quantities, but we do observe very strong contractivity for $\ME^s$, except in the case of $k$ small and $N$ large. In the latter case the problem is not very indefinite: and  GMRES iteration counts are modest even though the norm of $\ME^s$ is large (we give these for the case $N = 8$ in the column headed {\tt GMRES}). 
In most of the experiments in the checkerboard case, $\ME^s$ is contracting when $s$ is much smaller that $N^2$. {In Figure \ref{fig:norms}, we  plot  $\|\ME^s\|$  against  $s$ and observe that $\|\ME^s\| < 1$  for exponents  $s \ll N^2$.
}

\begin{figure}[t]
  \begin{center}
  \includegraphics[width=0.40\textwidth]{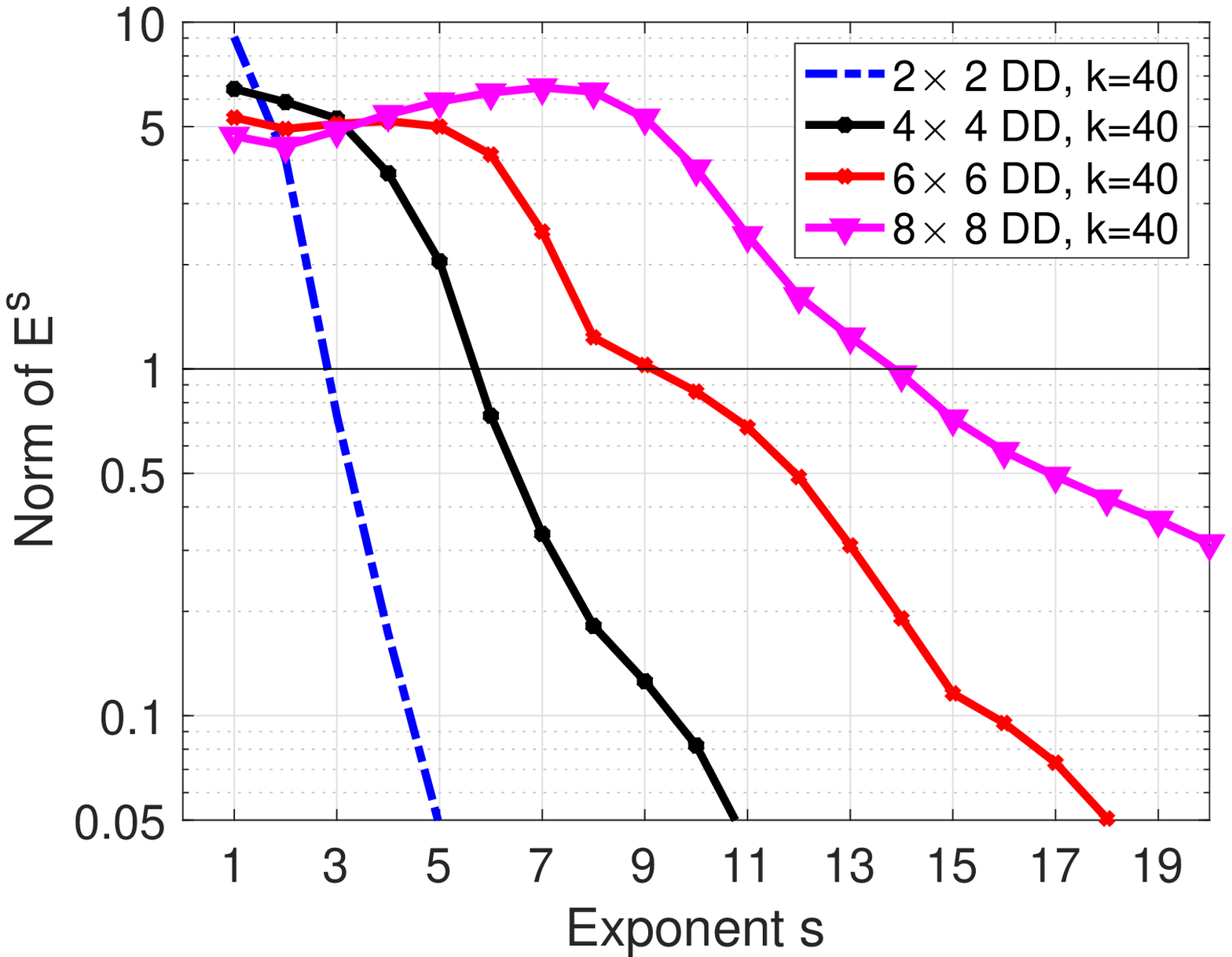}~~~~
    \includegraphics[width=0.40\textwidth]{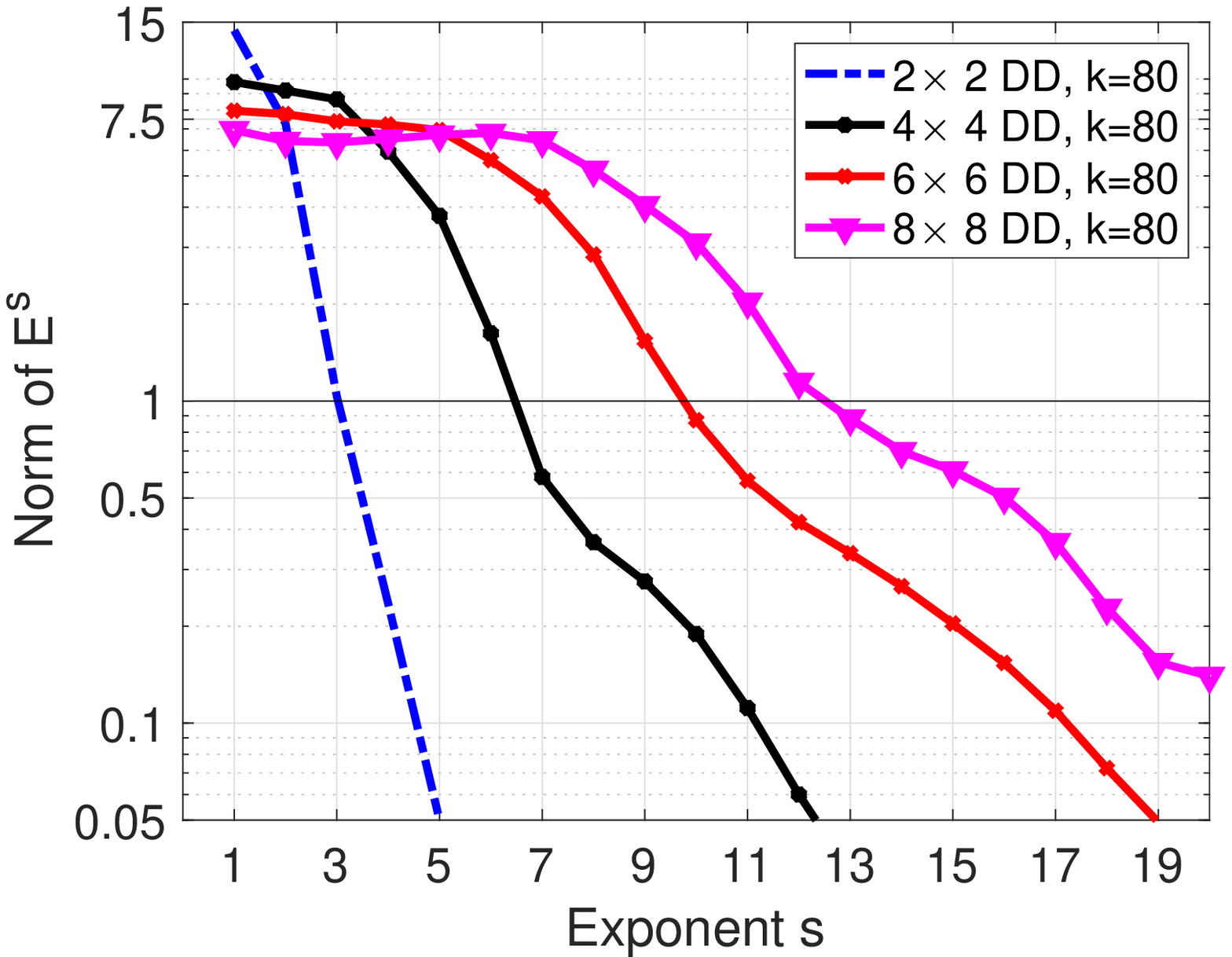}
\end{center}   
\caption{Norm of the power of the error propagation matrix (left: $k=40$, right: $k=80$)  \label{fig:norms}}
\end{figure}

\bigskip

\noindent {\bf Acknowledgement} SG thanks the {Section de Math\'{e}matiques}, University of Geneva for their hospitality during his visit in early 2020.
We gratefully acknowledge support from the UK Engineering and Physical Sciences Research Council Grants EP/R005591/1 (EAS) and   EP/S003975/1 (SG, IGG, and EAS). 

\vspace{-0.3cm} 

\end{document}